\documentclass[12pt]{article} 

\usepackage[letterpaper,textwidth=468pt,textheight=630pt,centering]{geometry}

\usepackage[utf8]{inputenc}
\usepackage[linesnumbered,ruled,vlined]{algorithm2e}
\usepackage{pdflscape}
\usepackage{comment,xcomment,amssymb,amsmath,color}
\usepackage{rotating,xypic,array}
\usepackage{latexsym,mathtools,amsthm}
\usepackage{subfig,rotating}
\usepackage{tikz,fancyvrb}
\usepackage{paralist}
\usepackage{booktabs}  
\usetikzlibrary{arrows.meta}
\usepackage{enumitem}
\usepackage{url} 
\usepackage{longtable}
\usepackage[section]{placeins}
\usepackage[longtable]{multirow}
\usepackage[colorlinks,citecolor=blue]{hyperref}
\usepackage{needspace}
\usepackage{etoolbox}

\AtBeginEnvironment{theorem}{\Needspace{8\baselineskip}}
\AtBeginEnvironment{lemma}{\Needspace{8\baselineskip}}
\AtBeginEnvironment{proposition}{\Needspace{8\baselineskip}}
\AtBeginEnvironment{corollary}{\Needspace{8\baselineskip}}
\makeatletter
\patchcmd{\@endtheorem}{\endtrivlist}
  {\@endparpenalty\@M\endtrivlist}{}{
    \PackageError{ispinor}{Cannot protect statement/proof page breaks}
      {Check the amsthm theorem-ending definition.}}
\AtBeginEnvironment{proof}{\@beginparpenalty\@M}
\makeatother

\newtheorem{theorem}{Theorem}[section]
\newtheorem*{theorem*}{Theorem}
\newtheorem{lemma}[theorem]{Lemma}

\newtheorem{corollary}[theorem]{Corollary}
\newtheorem*{corollary*}{Corollary}
\newtheorem{proposition}[theorem]{Proposition}
\theoremstyle{definition}
\newtheorem{definition}[theorem]{Definition}
\newtheorem{example}[theorem]{Example}

\theoremstyle{remark}
\newtheorem{remark}[theorem]{Remark}

\newcommand{\abs}[1]{\left\vert#1\right\vert}
\newcommand{\R}{\mathbb{R}}
\newcommand{\lie}[1]{\mathfrak{#1}}     

\newcommand{\Z}{\mathbb{Z}}

\newcommand{\C}{\mathbb{C}}

\newcommand{\CP}{\mathbb{CP}}

\newcommand{\hook}{\lrcorner\,}
\newcommand{\Pin}{\mathrm{Pin}}
\newcommand{\Spin}{\mathrm{Spin}}
\newcommand{\SU}{\mathrm{SU}}

\newcommand{\SO}{\mathrm{SO}}

\newcommand{\id}{\operatorname{Id}}   

\newcommand{\gl}{\lie{gl}}

\newcommand{\Span}[1]{\operatorname{Span}\left\{#1\right\}}

\newcommand{\eps}{\varepsilon}

\DeclareMathOperator{\Hom}{Hom}

\DeclareMathOperator{\Tr}{tr}
\DeclareMathOperator{\sign}{sign}

\DeclareMathOperator{\RT}{RT}

\newcommand{\scal}{\mathrm{scal}}

\newcolumntype{C}{>{$}c<{$}}
\newcolumntype{L}{>{$}l<{$}}
\newcolumntype{R}{>{$}r<{$}}

\usepackage{float} 
\usepackage{mathrsfs}
\newcommand{\odiv}{\mathbin{\oslash}}
\def\2{\mathrm{G}_2}
\def\8{\infty}
\def\0{\mathrm{o}}
\def\tq{\!\to\!}
\def\we{\wedge}
\def\La{\Lambda}
\def\RT{\mathrm{RT}}

\def\cd{\cdot}
\def\op{\oplus}
\def\O{\mathrm{O}}
\def\U{\mathrm{U}}

\def\sT{\mathscr{T}}
\def\sS{\mathscr{S}}
\def\sT{\mathscr{T}}
\def\sU{\mathscr{U}}

\def\FF{\mathbb{F}}
\def\ol{\overline}
\def\ge{\geqslant}
\def\le{\leqslant}
\def\vph{\vphantom{\int_p^d}}

\def\ba{\begin{array}}
\def\ea{\end{array}}

\def\QR{\mathrm{QR}}
\def\si{\sigma}
\def\Si{\Sigma}
\def\ti{\!\times}
\def\suml{\sum\limits}

\begin{document}
\parskip1pt

\title{Tournaments, flags, and spinors}

\author{Diego Conti and Simon Salamon}

\maketitle 
\begingroup
\renewcommand{\thefootnote}{}
\footnotetext{%
}
\endgroup

\abstract{It is observed that the dimension of the space of invariant
  spinors on the manifold $\SU(n)/T^{n-1}$ of full flags in $\C^n$
  equals the number $\RT(n)$ of labelled regular tournaments of order $n$, for
  $n$ odd. This space is decomposed under a double cover of the Weyl
  group, and realized as a pair of eigenspaces, first for the Dirac operator, then for a twisted version of it which acts equivariantly; this is used to show that multiplicities of irreducible
representations are even for $n\equiv 1\pmod 4$.}

\section*{Introduction}

This note pursues one of the problems addressed in
\cite{ArtachoSemmelmann}, which studied $G$-invariant
spinors on full flag manifolds $G/T$, where $G$ is a simple compact
Lie group and $T$ a maximal torus. Here, we are only concerned with
the case $G=\SU(n)$. The analogous problem for homogeneous spheres was studied in \cite{AgricolHofmannLawn}.

It was pointed out in \cite{Burstall_Salamon_1987} that
$\SU(n)$-invariant almost complex structures on the manifold
$\FF_n=\SU(n)/T^{n-1}$ correspond to tournaments with $n$ vertices, a
tournament being a digraph whose underlying graph is complete, so
$K_n$. In this setup, the almost complex structure is integrable if
and only if the tournament has no cycles. Based on a special case
described in \cite{Salamon_1985}, that setup was used in
\cite{Burstall_Salamon_1987} to help classify harmonic maps of $S^2$
into complex Grassmannians in which second fundamental forms of line
bundles over the domain define a quiver. That approach is in a sense
subsumed in Uhlenbeck's theory of uniton number \cite{Uhlenbeck}.

In the simplest example, $\SU(3)/T^2$ has a unique non-integrable
invariant almost complex structure up to sign, represented by a
3-cycle. This is familiar from the flag manifold's realization as the
nearly-K\"ahler twistor space of $\CP^2$. Our work highlights
generalizations of this almost complex structure, using the language
of graph theory.

A digraph is called \emph{balanced} if the in-degree (number of
incoming arrows) equals the out-degree (number of outgoing arrows) at
each vertex; its \emph{order} is the number of vertices. A balanced
tournament is also called a \emph{regular tournament}, and if its
order is $n$ both degrees equal $(n-1)/2$, and $n$ must be odd. There
is a total of $2^N$ labelled tournaments of order $n$, where
$N=\frac12n(n-1)$, ignoring the question of isomorphism.  Let $\RT(n)$
denote the number of the $2^N$ labelled tournaments that are
regular. As observed in \cite{McKay_1990}, $\RT(n)$ is the constant
term in
\begin{equation}\label{prod}
  \prod_{0\le j<k\le n-1}(x_j^{-1}x_k+x_k^{-1}x_j).
\end{equation}
The formula is used to tabulate values in Section \ref{sec:Ss}, and is
explained above Table \ref{Table:RT}.

Using the decomposition of the spinor representation under the inclusion
$\SU(N)\hookrightarrow \Spin(2N)$ defined by a regular tournament
$\sT$, we observe that a balanced subgraph of $\sT$ defines an
invariant spinor. But such a balanced subgraph $\sS$ defines a new
regular tournament $\sT\odiv\sS$ by reversing the arrows of $\sT$
along $\sS$. This is the basis of a simple but striking result of Akin
\cite[Corollary 4.5]{Akin_2020}, asserting that the number of balanced
subgraphs of $\sT$ equals $\RT(n)$. After filling in details in
Section~\ref{sec:acs} and reinterpreting Akin's result in our setting (Theorem~\ref{akin}), we deduce Corollary~\ref{cor:dim}, effectively an observation, which states that the space of invariant spinors on
$\SU(n)/T^{n-1}$ has dimension $\RT(n)$.

We discuss the way in which the choice of a regular tournament $\sT$,
or `regular' almost complex structure, imposes an algebra
structure on the space of invariant spinors. This provides refined
counts of spinors by degree that can be used to distinguish
isomorphism classes of regular tournaments for $n\ge7$. This is
explained in Section \ref{sec:Ss}, and illustrated in Section
\ref{sec:low}, which includes an analysis of the cases $n=5$ and
$n=7$.

Section~\ref{sec:weyl} introduces the action of a double cover of the
Weyl group on the space of invariant spinors and computes its
character; Section~\ref{sec:decomposition} deduces the decomposition
into irreducible spin representations, and Section~\ref{sec:dirac}
identifies the space of invariant spinors with the sum of two
eigenspaces of the Dirac operator. A similar decomposition for a twisted version of the Dirac operator is used to show that self-associate representations, and all representations for $n\equiv 1\pmod 4$, appear with even multiplicity.

\section{Invariant almost complex structures}\label{sec:acs}

Notwithstanding the general theory of invariant spinors on flag
manifolds $G/T$, described in \cite{ArtachoSemmelmann} using weights,
we are only concerned with the case $G=\SU(n)$ and the flag manifold
\[ \FF_n=\SU(n)/T^{n-1}=\U(n)/T^n,\]
which has real dimension $2N$, where $N={n\choose2}$.  We fix an
$\SU(n)$-invariant Riemannian metric on $\FF_n$, equivalently an inner
product on $\R^{2N}$ invariant by the action of $T^n$, though the
choice does not affect the outcome.

Our starting point is that a manifold with an orthogonal almost
complex structure (that is, a $\U(N)$-structure) is spin if and only
if the canonical bundle $\kappa\cong\Lambda^{N,0}$ has a square
root. In this case, one can express the associated spinor representation
in terms of $(0,q)$-forms, as in Hitchin's influential paper
\cite{Hitchin_1974}, by
\begin{equation}
 \label{eqn:hitchinspinors}
 \Delta = \kappa^{1/2}\otimes(\La^{0,0}+\La^{0,1}+\dotsc + \La^{0,N}).
\end{equation}
Here, $+$ denotes an orthogonal sum.
Recall that $\lie{u}(n)\cong[\,\ol{\C^n}\otimes\C^n]$, in which single
square brackets indicate `the real vector space whose complexification
is what is inside' \cite{Salamon_1989}. Under $\lie t^n$, we have the
decomposition
\[\C^n=L_0+\cdots+ L_{n-1}\]
where $+$ denotes direct sum, so
\[\begin{array}{rcl} \lie u(n)
&=& [(\ol L_0+ \dotsc +\ol L_{n-1})\otimes(L_0+ \dotsc + L_{n-1})]\\[5pt]
&=& [\sum\limits_i \ol{L_i}\otimes L_i] + [\sum\limits_{i\neq j}\ol{L_i}\otimes L_j]\\[5pt]
  &=& \lie t^n+ \sum\limits_{i<j}[\ol{L_i}\otimes L_j+\ol{L_j}\otimes L_i]\\[8pt]
  &=& \lie t^n+ \sum\limits_{i<j} [\![\ol{L_i}\otimes L_j]\!].
\end{array}\]
Here, $[\![\ol{L_i}\otimes L_j]\!]$ denotes the real 2-dimensional vector space underlying
\[\ol{L_i}\otimes L_j\cong\Hom(L_i,L_j).\]
It will be convenient to identify $\lie u(n)/\lie t^n$ with a generic
\emph{co}tangent space of the flag manifold $\FF_n$.  A $\lie
t^n$-equivariant almost complex structure $J$ must leave each of the
2-dimensional spaces invariant. Moreover, the space of $(1,0)$-forms
for $J$ on $[\![\ol L_i\otimes L_j]\!]$ is a $\lie t^n$-equivariant
complex line inside
\[\ol{L_i}\otimes L_j +\ol{L_j}\otimes L_i\]
so it must be one summand or the other. Indeed, $J$ amounts to a
choice of orientation of $[\![\ol L_i\otimes L_j]\!]$.

An invariant almost complex structure $J$ can therefore be identified
with a directed graph with $n$ nodes, where $i\to j$ if and only if
$\ol{L_i}\otimes L_j$ has type $(1,0)$.

The discussion above implies

\begin{lemma}
\label{lemma:invariantalmostcomplexstructures}
Invariant almost complex structures on $\FF_n$ are in one-to-one
correspondence with tournaments of order $n$.
\end{lemma}

A \emph{cycle} in a digraph is a connected balanced subgraph for which the
in-degree (and so out-degree) of each vertex appearing in the subgraph
is exactly one. (Some authors refer to this as a \emph{simple cycle}.)
The following is completely standard, but worth keeping in mind:

\begin{lemma}\label{std}
A finite digraph is balanced if and only if its set of arrows can be
partitioned into cycles.
\end{lemma}

\begin{proof}
By induction on the number of arrows. If there are no arrows, there is
nothing to prove. Otherwise, let $e_1$ be an arrow from $i_0$ to
$i_1$. Since the digraph is balanced, there is some arrow going out of
$i_1$, say $e_2$. This gives an infinite sequence $\{e_j\}$; because
the digraph is finite, at some point a vertex must be repeated, and
this gives a cycle.  Removing the cycle, one obtains another
balanced digraph, to which the induction hypothesis can be applied.
\end{proof}

To establish notation, we shall describe some tournaments. We shall label the vertices of a tournament of order $n$ using the
integers from the set $\{0,1,2,\ldots,n-1\}$. Here then is an example:

\begin{definition}\label{circ}
Let $n\ge3$ be odd. The \emph{circulant tournament}, denoted in this
paper by $\sT_n$, is the digraph in which $i\to j$ is an edge if and
only if $j-i\equiv 1,\ldots, (n-1)/2\pmod n$.
\end{definition}

If $n$ is prime, then $\sT_n$ can be written as the union of
$(n-1)/2$ $n$-cycles. Here is a description based on  Lemma \ref{std}:
\begin{equation}\label{3579}\ba{rcl}
\sT_3 &=& 012\\
\sT_5 &=& 01234\cd02413\\
\sT_7 &=& 0123456\cd0246135\cd0362514\\
\sT_9 &=& 012345678\cd024681357\cd036\cd147\cd258\cd 048372615
\ea
\end{equation}
Dots merely indicate unions of cycles, whose parentheses we omit.
As an example, the notation $02413$ indicates the 5-cycle
\[ 0\to2\to4\to1\to3\to0,\]
which can also be written $24130$ etc. Note that we do not repeat the
choice of initial vertex in the abbreviated notation.

It is important to realize that there are in general many valid
decompositions of a tournament into disjoint cycles. For example,
\[ \sT_5=0123\cd024\cd134\]
reflecting the partition $10=4+3+3$. The circulant tournament $\sT_n$
is very easy to work with because one can construct cycles by adding
$1,2,\cdots,\frac12(n-1)$ to the labels, without reference to a
diagram or crib sheet. We shall introduce inequivalent tournaments in
the next section.  \medbreak

The next result is central to our approach.

\begin{lemma}\label{tacs}
The tournament associated to an invariant almost complex structure on
$\FF_n$ is regular if and only if $\kappa$ is trivial.
\end{lemma}

\begin{proof}
The space $\La^{1,0}$ of $(1,0)$-forms is a direct sum of the
lines $\ol{L_i}\otimes L_j$, so there is an isomorphism
\[\kappa\ \cong\ \bigotimes_{i\to j} \ol{L_i}\otimes L_j.\]
If $\sT$ is regular, each occurrence of $L_i$ on the right-hand side
will be matched by an occurrence of $\ol{L_i}$ (and their
multiplicities will also be equal). Thus, the product is trivial.

Conversely, suppose that $\kappa$ is trivial. Treat $\ol{L_i}$ on the
right-hand side as $L_i^{-1}$. Each arrow has a head and a tail, so
the total number of negative powers must equal the total number of
positive powers. Now tensor the right hand side by a power of the
trivial tensor product
\[ \lambda=L_0\otimes\cdots\otimes L_{n-1}\cong\La^n(\C^n)\]
so that only non-negative powers of all $L_i$ persist. Since the
product remains trivial, it must itself be a power of $\lambda$, as
such a power is the only trivial 1-dimensional representation of
$T^{n-1}$. It follows that the original product was (up to permutation
of factors) $\ol\lambda{}^{(n-1)/2}\otimes\lambda^{(n-1)/2}$, and the in-
and out-degrees are equal at every node.
\end{proof}

We will say an invariant almost complex structure is \emph{regular} if
$\kappa$ is trivial. It follows that (for any $n$ odd), $\sT_n$
determines a regular invariant almost complex structure on $\FF_n$.
This will enable us to count the dimension of spaces of invariant
spinors.\smallskip

The following result is well known, but we can prove it pursuant to
Lemma~\ref{tacs}.

\begin{proposition}\label{spin}
For any $n$, the full flag manifold $\FF_n$ is spin.
\end{proposition}

\begin{proof}
  If $n$ is odd, $\FF_n$ admits an invariant orthogonal almost complex
  structure with $\kappa$ trivial, equivalently an $\SU(N)$
  structure. The latter lifts to an inclusion $\SU(N)\to \Spin(2N)$.

If $n$ is even, we need to show that for some invariant almost complex
structure $\kappa$ admits a square root. Label the nodes
$0,1,2,\ldots,n-1$, and choose a regular tournament on the complete
subgraph $K_{n-1}$ generated by $1,\ldots, n-1$. Then add all the
directed edges $0\to k$ for $1\le k\le n-1$, so as to define a
so-called reducible tournament. In the regular tournament on
$\{1,\dots,n-1\}$ each vertex has equal in- and out-degree, so its
contribution $\bigotimes_k(\ol{L_k}\otimes L_k)^{(n-2)/2}$ to the
canonical bundle is trivial, and
  \[\kappa \cong \ol{L_0}^{\,n-1}\otimes L_1\otimes\cdots\otimes L_{n-1}\cong
  \ol{L_0}^{\,n-1}\otimes\ol{L_0}=\ol{L_0}^{\,n},\]
in which the middle isomorphism arises from the fact that
$L_0\otimes\cdots\otimes L_{n-1}=\La^n(\C^n)$ is a trivial line
bundle. Thus $\kappa^{1/2}\cong \ol{L_0}^{\,n/2}$, which exists
because $n$ is even.
\end{proof}

\section{Spinors and subgraphs}\label{sec:Ss}
In light of \eqref{eqn:hitchinspinors}, if we identify an almost
complex structure for which $\kappa$ is trivial, thereby reducing the
structure at each point to $\SU(N)$, we obtain
\[\Delta = \La^{0,0}+\dotsc + \La^{0,N}.\]
One can also replace $(0,q)$ forms by ones of type $(p,0)$, since
$\La^{0,q}\cong\La^{0,q}\otimes\La^{N,0}\cong\La^{N-q,0}$.
Explicitly, Clifford multiplication takes the form of an equivariant map
\[(\La^{1,0}+\La^{0,1})\times \Delta\to \Delta,\]
induced by
\[\begin{array}{l}
\La^{1,0}\otimes \La^{0,q}\to \La^{0,q-1}, \quad (\alpha,\beta)\mapsto -\sqrt2\alpha\hook\beta\\[3pt]
\La^{0,1}\otimes \La^{0,q}\to\La^{0,q+1}, \quad (\alpha,\beta)\mapsto \sqrt2\,\alpha\we\beta.
\end{array}\]
The signs are chosen so that a real tangent vector $v$ satisfies
$v\cdot v=-\abs{v}^2$.

An invariant spinor on $\FF_n$ is a section of the associated vector
bundle that is invariant by the action lifted from that of $\SU(n)$ on
$\FF_n$. It will be covariant constant for the canonical connection on
the spin bundle, but not for any Levi-Civita connection. There are no
nonzero invariant spinors if $n$ is even, as can be reproved with our technique, see Lemma~\ref{lemma:neven} below.

For each arrow $i\to j$ in $\sT_n$, let $\alpha_{ij}$ be the generator of $\overline L_i\otimes L_j\cong\Hom(L_i,L_j)$ corresponding to the elementary matrix with $1$ at the entry $(j,i)$. Then $\alpha_{ij}$ has type $(1,0)$ by construction. Given a subgraph $\sU$ in $\sT_n$, we define
\[\psi_\sU=\prod_{i\to j\text{ in }\sU} \overline\alpha_{ij}.\]
This gives a decomposable element in $\La^{0,q}$, with $q$ being the number of edges in $\sU$. Since the sign of the product depends on the order of the factors, we will assume that an order on the edges in the tournament has been fixed, though the particular choice is irrelevant.

\begin{lemma}[\cite{ArtachoSemmelmann}]
\label{lemma:neven}
If $\FF_n$ has nonzero invariant spinors, then $n$ is odd.
\end{lemma}
\begin{proof}
Fix an invariant almost complex structure associated to a tournament $\sT$; a spinor is a section of $\kappa^{1/2}\otimes\La^{0,\bullet}$. As a representation of $\lie t^n$, $\La^{0,\bullet}$ decomposes into
the direct sum of $1$-dimensional representations, each corresponding to a decomposable form in some $\La^{0,q}$. Accordingly, an invariant spinor is a linear combination of invariant elements of the form $\epsilon\otimes\psi_\sU$, where $\epsilon$ is a generator of $\kappa^{1/2}$ and $\sU$ varies among subgraphs of $\sT$, balanced or not. The square of $\epsilon\otimes\psi_\sU$ effectively belongs to
\[\bigotimes_{i\to j \text{ in } \sT\setminus\sU} \overline L_i\otimes L_j \otimes \bigotimes_{i\to j\text{ in } \sU}\overline{L_j}\otimes L_i,\]
which is the canonical bundle of the almost complex structure obtained by reversing the arrows of $\sT$ contained in $\sU$. Arguing as in Lemma~\ref{tacs}, the fact that this canonical bundle is trivial implies that the in-degree and out-degree coincide at each vertex, which is only possible if $n$ is odd.
\end{proof}
In light of Lemma~\ref{lemma:neven}, from now on we will assume that $n$ is odd. Having fixed the regular almost complex structure defined by $\sT_n$, we identify the space $\Xi_n$ of invariant spinors as $\lie t^{n-1}$-invariant elements of $\La^{0,\bullet}$. The inclusion in  $\La^{0,\bullet}$ makes $\Xi_n$ into an algebra, in a way that depends on the fixed almost complex structure used for the identification.

\begin{lemma}
\label{lemma:basis}
If $\sU$ is a balanced subgraph, $\psi_\sU$ is invariant; the $\psi_\sU$ thus obtained form a basis of the space of invariant spinors.
\end{lemma}
\begin{proof}
The argument of Lemma~\ref{tacs} shows that invariance of $\psi_\sU$ is equivalent to $\sU$ being balanced; that this gives a basis follows from the fact that each $\psi_\sU$ spans  a $1$-dimensional representation of $\La^{0,\bullet}$.
\end{proof}

Recall from \cite{Chevalley_1996} that a spinor is \emph{pure} if the subspace of the complexified tangent space that annihilates it is maximally isotropic, hence the space of vectors of type $(1,0)$ relative to some almost complex structure; this defines a bijection between orthogonal almost complex structures and pure spinors up to scale.
\begin{lemma}
\label{lemma:pure}
If $\sU$ is a balanced subgraph, $\psi_\sU$ is invariant and pure; the almost complex structure defined by $\psi_\sU$ is regular, and corresponds to the tournament  $\sT_n\odiv \sU$ obtained from $\sT_n$ by reversing the arrows of $\sU$.
\end{lemma}

\begin{proof}
Given $\sU$, we compute the space
\[V=\{v\mid v\cdot\psi_\sU=0\}.\]
It is clear that
\begin{equation*}
\begin{split}
\overline\alpha_{ij}\cdot\psi_\sU&=\sqrt2\overline\alpha_{ij}\we \psi_\sU=0 \text{ if $i\to j$ in $\sU$}\\
 \alpha_{ij}\cdot\psi_\sU&=-\sqrt2\alpha_{ij}\hook \psi_\sU=0 \text{ if $i\to j$ is not in $\sU$}
 \end{split}
\end{equation*}
Thus, $V$ has maximal dimension and $\psi_\sU$ is pure. Its space of
$(1,0)$-vectors is obtained from the space of $(1,0)$-vectors for the
almost complex structure defined by $\sT_n$ by conjugating the
components corresponding to elements of $\sU$. In particular the
canonical bundle is generated by
\begin{equation}\label{wedge}
  \prod_{i\to j\in\sU} \overline\alpha_{ij}\;\;\hbox{\large$\we$}\kern-8pt
  \prod_{i\to j\in\sT_n\setminus \sU} \alpha_{ij}.
\end{equation}
Arguing as in Lemma~\ref{tacs}, we see that both factors are
invariant, hence $\kappa$ is invariant and the almost complex
structure defined by $\psi_\sU$ is regular. Thus, it is the almost
complex structure defined by a regular tournament, which by construction is $\sT_n\odiv \sU$.
\end{proof}

\begin{lemma}
\label{lemma:trivialcanonicalbundle}
The almost complex structure defined by a pure invariant spinor is regular.
\end{lemma}

\begin{proof}
By \cite[III.3.2]{Chevalley_1996}, we have a squaring map
\[\Delta \otimes\Delta\to\La^{\bullet}\otimes \La^{0,N},\]
with $\La^\bullet$ denoting the complexified exterior algebra; if $V$ is the maximal isotropic space determined by a spinor $\psi$,
the image of $\psi\otimes \psi$ lies in $\det
V\otimes\La^{0,N}$. Since $\psi\otimes\psi$ is invariant and
$\La^{0,N}=\ol\kappa$ is trivial (the reference structure $\sT_n$
being regular), this implies that $\det V$ is also trivial, i.e. the
almost complex structure defined by $\psi$ has trivial canonical
bundle.
\end{proof}

Our point of view leads to an alternative proof of a result of Akin,
which enables one to count invariant spinors on $\FF_n$. The result
itself is remarkable because it asserts that to list \emph{all}
regular tournaments of order $n$, one only need look inside any
one. Expressed loosely, a regular tournament is a perfect or universal
object.

\begin{theorem}[\protect{\cite[Corollary 4.5]{Akin_2020}}]\label{akin}

The number of balanced subgraphs of a regular tournament of order $n$
is precisely $\RT(n)$.
\end{theorem}

\begin{proof}
Fix a regular tournament $\sT$. Lemma~\ref{lemma:pure}, Lemma~\ref{lemma:trivialcanonicalbundle} and
Lemmas~\ref{lemma:invariantalmostcomplexstructures}, \ref{tacs} (in
this order) determine maps $f,g,h$ that fit in the following
commutative diagram:
\[\xymatrix{\{\text{balanced subgraphs of } \sT\}\ar[r]\ar[d]^f & \{\text{regular tournaments on the $n$ nodes}\}\ar[d]^{h}\\
 \{\text{pure invariant spinors}\}/\C^* \ar[r]^g&
  \left\{\parbox{3cm}{\centering regular invariant\\ almost complex structures}\right\}
 }
\]
In this diagram, $f$ is injective and $h$ is a bijection; the map $g$ is a special case of the correspondence between almost complex structures and pure spinors, hence also a bijection. Furthermore, the argument leading to Lemma~\ref{lemma:invariantalmostcomplexstructures} shows that every invariant almost complex structure is obtained from the fixed one by reversing the arrows of a subgraph; by Lemma~\ref{tacs} the resulting structure is regular precisely when the subgraph is balanced. Thus $g\circ f$ is surjective, and all the arrows in the diagram are bijections.
\end{proof}

Lemma~\ref{lemma:basis} immediately implies

\begin{corollary}\label{cor:dim}
The dimension of the space of invariant spinors $\Xi_n$ on $\FF_n$ is $\RT(n)$.
\end{corollary}

Although there is no closed formula for $\RT(n)$, one can tabulate
values by exploiting \eqref{prod}, in which a term $x_j^{-1}x_k$
signals the presence of the arrow $j\to k$. Each monomial evaluating
to $1$ is a product of such `arrows' in which exactly one of $j\to k$
and $k\to j$ appears for each pair $j<k$. This specifies a
tournament and vice versa.

\begin{table}[ht]
\small
\centering
\begin{tabular}{r|p{0.7\textwidth}}
$n$ & $\RT(n)$ \\
\hline
1 & 1 \\
3 & 2 \\
5 & 24 \\
7 & 2 640 \\
9 & 3 230 080 \\
11 & 48 251 508 480\\
13 & 9 307 700 611 292 160\\
15 & 24 061 983 498 249 428 379 648\\
17 & 855 847 205 541 481 495 117 975 879 680\\
19 &  427 102 683 126 284 520 201 657 800 159 366 676 480\\
21 &3 035 991 776 725 501 434 069 099 002 640 396 043 332 019 814 400
\end{tabular}
\normalsize
\caption{Counts of labelled regular tournaments\label{Table:RT}}
\end{table}

Exact values of $\RT(n)$ for $n\le21$ were found by  \cite{mckay1983applications} summing over
roots of unity, and extended to $n\le37$
by means of a recurrence on out-degree sequences \cite[Appendix]{Isaev_McKay_Zhang_2024}.
There is an asymptotic expression
\[ f(n)=\sqrt{n/e}\left(\frac{2^{n + 1}}{\pi n}\right)^{(n - 1)/2},\]
with multiplicative error \[1+O(n^{-1/2+\eps})\] for any $\eps>0$
\cite{McKay_1990}. For example, the value $f(21)\approx 2.998\times
10^{48}$ can be compared with the final row of Table \ref{Table:RT}.

Incidentally, the high connectivity in the unitary case is reflected
by the value $\RT(17)\approx 8.56\times 10^{29}$ corresponding to
dimension $N=136$. This already exceeds a predicted dimension of
$2\times10^{29}$ (to one significant figure) for the space of
invariant spinors on $E_8/T^8$, which has dimension $240$. This raises
the question of whether there is a graph theoretic version that is
applicable to other full flag manifolds, perhaps by combining the Lie
combinatorics of \cite{ArtachoSemmelmann} and \cite{MoroianoSchwann}.

We stress the fact that $\RT(n)$ counts all regular tournaments, without
identifying `equivalent' tournaments.  Formally, an \emph{isomorphism}
between two digraphs is a bijection of the sets of nodes that induces
a bijection of arrows; an automorphism is an isomorphism of a digraph
to itself. Counts of isomorphism classes of tournaments appear in
\cite{Moon_1968,Chamberland_Herman_2015}. Regular tournaments taken up
to automorphisms can be listed using the programme {\tt nauty}, see
\cite{McKay_Piperno_2014}.

\smallskip
Having fixed a regular tournament $\sT$ of order $n$ and associated
$\SU(N)$ structure on $\FF_n$, we obtain a decomposition
\[ \Xi_n=\bigoplus_{q=0}^N \Xi_{n,q},\]
with the notation $\Xi_{n,q}=\Xi_n\cap\La^{0,q}$. For example,
$\dim\Xi_{n,0}=1$, corresponding to $\La^{0,0}$ and the empty
subgraph. Similarly $\dim\Xi_{n,N}=1$, corresponding to $\sT$ itself
and the trivial space $\La^{0,N}=\ol\kappa$. More generally,
`Poincar\'e duality' holds:

\begin{lemma}\label{Pd}
  $\dim\Xi_{n,q}=\dim\Xi_{n,N-q}$.
  \end{lemma}

\begin{proof} If $\sS$ is a balanced subgraph of $\sT$,
    so is its complement $\sT\setminus\sS$, and this establishes an
    isomorphism $*\colon\Xi_{n,q}\to\Xi_{n,N-q}$.
\end{proof}

\section{Regular tournaments of low order}\label{sec:low}

This section is designed to give a taste of the classification of
regular tournaments of order at most $9$, and to illustrate some of
the preceding theory and applications. Understanding its concrete
examples was also a determining factor in leading the authors to
develop the link with spinors. We have seen that each regular
tournament determines an almost complex structure with trivial
canonical bundle, and its balanced subgraphs give a basis for the
space of invariant spinors on $\FF_n$, whose elements
can be listed by degree.

As a starting (and reference) point, we shall use the circulant
tournament $\sT_n$ introduced in Definition~\ref{circ}. Its vertices will
again be labelled $0,1,2,\ldots,n-1$.\smallskip

There is little to say in the case $n=3$. Indeed, $\sT_3$ is
represented by the 3-cycle $012$, or triangle
\[ 0\to 1\to 2\to 0\]
which defines a non-integrable almost complex structure $J$ on $\FF_3$
mentioned in the Introduction. There are two balanced subgraphs, the
empty one and $\sT_3$ itself, corresponding to $J$ and $-J$. Up to
obvious isomorphisms, there are exactly two tournaments of order 3,
balanced (non-integrable) and cycle-free (integrable).

The space of $\SU(3)$-invariant spinors on $\FF_3$ has dimension 2. It
is generated by forms of degree $0$ and $3$, and contains two
\emph{Killing} spinors for the Levi-Civita connection of the
nearly-K\"ahler metric.

\subsubsection*{The regular tournament of order 5}\label{ss5}

  The fact that $\RT(5)=24$ was known to the authors of
\cite{ArtachoSemmelmann}. Assuming this result for the moment, it
follows that any regular tournament of order $5$ is isomorphic to
\[\sT_5=01234\cd02413.\]
Indeed, when the symmetric group $\Si_5$ acts on $\sT$, its stabilizer
is the subgroup $\Z_5$ generated by $01234$ (whose square is its
`cofactor' $02413$), so the orbit has size $5!/5=24$. It is a
coincidence that $24$ is the dimension of the isometry group $\SU(5)$,
which by construction acts trivially on the elements of $\Xi_5$.

Although one can list all balanced subgraphs of $\sT_5$ by hand, this
exercise is already involved, especially when one does not know that
their total is 24. However, the formula \eqref{prod} enables one to
construct (or, in the authors' case, verify) Table \ref{deg5}, which
lists the 24 balanced subgraphs of $\sT_5$.

{\small\centering
  \begin{table}[H]
\[\begin{array}{clr}
\hbox{deg} & \hbox{cycle structure}& \hbox{\#}\\\hline
0 & \varnothing & 1\\
3 & 013,\ 023,\ 024,\ 124,\ 134 & 5\\
4 & 0123,\ 0124,\ 0134,\ 0234,\ 1234 & 5\\
5 & 01234,\ 02413 & 2\\
6 & 013\cdot024,\ 013\cdot124,\ 023\cdot124,\ 023\cdot134,\ 024\cdot134 & 5\\
7 & 0123\cdot024=0124\cdot023,\ 0134\cdot023=0234\cdot013 &\\[-3pt]
  & 0234\cdot124=1234\cdot024,\ 0123\cdot134=1234\cdot013 &\\[-3pt]
  & 0124\cdot134=0134\cdot124 & 5\\
10 & 01234\cdot02413 & 1\\\hline
& & 24
\end{array}\]
\vspace{-10pt}
\caption{Counts of balanced subgraphs of $\sT_5$ by length}\label{deg5}
  \end{table}
}

One can visualize $\sT_5$ as a directed regular pentagon $01234$
together with an inscribed pentagram $02413$. Each balanced subgraph
with 7 edges can be decomposed in two ways as the union of a 4-cycle
and a 3-cycle. For example, both $0123\cd024$ and $0124\cd023$
determine the subgraph $\sS$ with edge set
\begin{equation}\label{sub7}
  0\tq1,\;\;\;0\tq2,\;\;1\tq2,\;\;2\tq3,\;\;2\tq4,\;\;3\tq0,\;\;4\tq0.
\end{equation}
The table confirms that there is indeed a total of 24 labelled regular
tournaments of order $5$, each of which is formed by reversing arrows
along $\sS$. For example, the tournament associated to \eqref{sub7} is
\[ \sT_5\odiv\sS = 03421\cd  04132.\]
Writing  $01$ for a generator of $\overline {L_0}\otimes L_1$ (denoted $\alpha_{01}$ in Section~\ref{sec:Ss}), and so on, the canonical bundle of the almost complex structure defined by $\sT_5$ is generated by
\[ \pm\;(01\we12\we23\we34\we40)\we
 (02\we24\we41\we13\we30);\]
the sign ambiguity arises because any two factors in the exterior
product can be swapped.

The correspondence of Lemma~\ref{lemma:basis} associates to $\sS$ a
spinor which we can write as
\[ \psi_{\sS}= \pm\>01\we02\we12\we23\we24\we30\we40\in \La^{0,7};\]
this time, $01$ represents a $(0,1)$-form
($\overline\alpha_{01}=-\alpha_{10}$ in the language of
Section~\ref{sec:Ss}), which has type $(1,0)$ in the almost complex
structure defined by $\sT\odiv\sS$.

Since any regular tournament of order $5$ is isomorphic to 
$\sT_5$, the decomposition
\[ \RT(5)=24=1+5+5+2+5+5+1\]
according to the length $q$ of balanced subgraphs (equivalently, the
degree of decomposable forms) is independent of the choice of
tournament. The balanced subgraphs are generated by 12 cycles of
lengths $3,4,5$. Put another way, the associated subalgebra of $\Lambda^{0,\bullet}$ has
$12=5+5+2$ generators, of respective degrees $3,4,5$.

As we shall demonstrate, uniqueness by order fails when
$n\ge7$. However, the counts of balanced subgraphs are always
consistent with Lemma~\ref{Pd}, and the number of 3-cycles is always
equal to $n(n^2-1)/24$, see \cite{Moon_1968}.

\subsubsection*{Regular tournaments of order 7}\label{ss7}

The number of regular 7-tournaments is $\RT(7)=2640$, divided (we
shall see) into three isomorphism classes.

\begin{figure}[!ht]
    \centering
      \includegraphics[width=.46\textwidth]{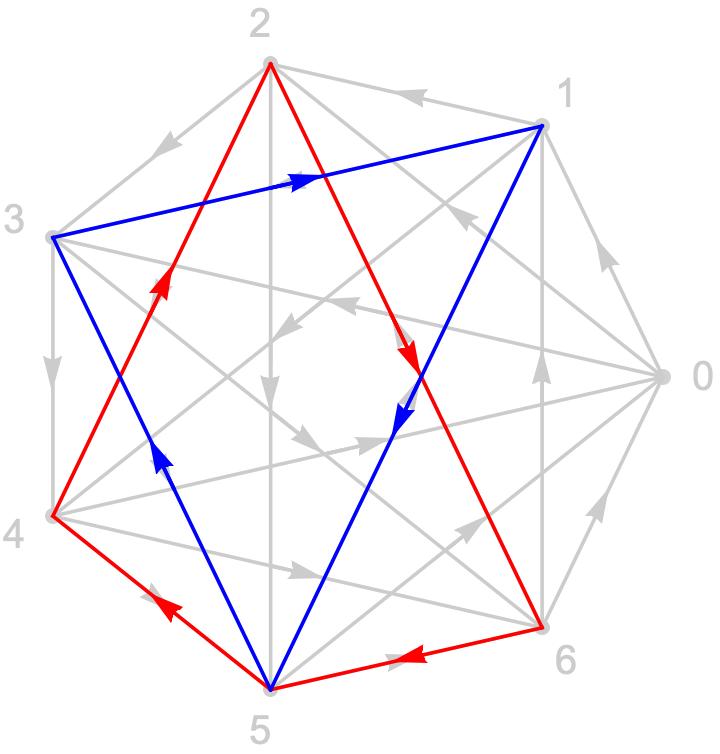}
    \caption{A modified 7-tournament\label{fig:septa}}
  \end{figure}

Our reference class is represented by
\begin{equation}\label{ref}
  \sT_7=0123456\cd0246135\cd0362514
\end{equation}
as in \eqref{3579}, observing that this is fixed by each of the three
7-cycles in the decomposition \eqref{ref}. A second tournament is
\[\ba{rcl} \sT_7'=\sT_7\odiv\sS
&=& 014\cd234\cd025\cd315\cd036\cd126\cd465\\
&=& 0126534\cd0231546\cd0361425.\ea\]
It can be obtained by reversing the arrows of $\sT_7$ on the edges of
its balanced subgraph
\[ \sS=2456\cdot 135\]
(chosen more or less at random) which is depicted in Figure~\ref{fig:septa}.

Although $\sS$ happens to be connected, any regular tournament of order
at least $7$ contains balanced subgraphs that are disconnected (for
example, consisting of two 3-cycles).

Camion's theorem \cite{Camion} states that a strongly connected
tournament (meaning there is a directed path between any two vertices)
possesses a Hamiltonian cycle. Kelly's conjecture states that every
regular tournament is the union of $(n-1)/2$ Hamiltonian cycles. It
has been verified for $n$ up to $9$ \cite{BangJensenGutin}, and is
known to be true for sufficiently large $n$ \cite{KuhnOsthus_2013}.

In fact, $\sT_7'$ is \emph{isomorphic} to the tournament
\[ \sT_7\odiv0362514=0123456\cdot 0246135\cdot 0415263.\]
This is the so-called quadratic residue tournament $\sT_7^\QR$
(sometimes called Paley tournament) that can be defined for $n$ a prime
with $n\equiv3$ mod $4$ by the rule
\[i\to j \iff  j-i\text{ is a square
modulo $n$}.\] This is well defined because $-1$ is not a square.

The authors call $\sT_7'$ the `$\2$ tournament', since its
decomposition into 3-cycles matches the 3-form
\[e^{014}+e^{234}+e^{025}+e^{315}+e^{036}+e^{126} - e^{456}
= \omega^1\we e^4+\omega^2\we e^5+\omega^3\we e^6 - e^{456}\] whose
stabilizer is conjugate to the subgroup $\2$ of $\SO(7)$. (Here
$\{\omega^i\}$ is a rescaled orthonormal basis of the subspace $\La^2_+$ of
self-dual 2-forms, and reflects the decomposition
$\R^7=\La^1\op\La^2_+$ relative to the subgroup $\SO(4)$ of $\2$
\cite{Salamon_1989}.) The link between these objects
merely provides a mnemonic for distinguishing the tournaments up to
isomorphism.

The third isomorphism class is represented by
\[ \sT_7''=\sT_7' \odiv 465 = 0123645\cd0261534\cd0314256\]
and can be remembered via the 3-form whose stabilizer is split $\2$!
In contrast to $\sT_7'$ and $\sT_7''$, the circulant tournament
$\sT_7$ cannot be expressed as a union of edge-disjoint $3$-cycles.
This is because $3$-cycles in the circulant tournament are translates of $014$ or $024$, hence can have either zero or two edges in
common with the cycle $0246135$; however, the latter has an odd number of
edges.

Table~\ref{deg7} lists the numbers of balanced subgraphs of
7-tournaments by length (equivalently, dimensions of invariant spinors
by degree) for the three isomorphism classes. In view of
Lemma~\ref{Pd}, there is no need to exceed $[N/2]$. The sum of the
entries of each row equals $1320=\frac12\RT(7)$.

\def\hq#1{\hbox to 24pt{\hfill$#1$}\phantom.} 
  \begin{table}[H]
\[\begin{array}{c|rrrrrrrrr}
q\tq &   \hq0&\hq3&\hq4&\hq5&\hq6&\hq7&\hq8&\hq9&\hq{10}\\\hline
\sT_7 & 1&14&28&28&^{28}91&^{17}220&245&266&427\vph\\
\sT_7' & 1&14&21&42&^{35}105&^{24}171&252&336&378\\
\sT_7'' &   1&14&24&36&^{35}102&^{15}183&249&330&381\\
\end{array}\]
\vspace{-5pt}
\caption{Counts of balanced subgraphs of $\sT_7,\sT_7',\sT_7''$ by length}\label{deg7}
\end{table}

By Lemma~\ref{std}, any balanced subgraph is generated by simple
cycles, and these are necessarily of length at most $n$. Some of these
subgraphs (for $q=6,7$ when $n=7$) will themselves be products of
smaller cycles. Superscripts in the entries count the numbers of
(simple) cycles. For example, for $\sT_7$, there will be $28$ 6-cycles
(and $91-28=63$ products of two 3-cycles), and $17$ 7-cycles (and
$220-17=203$ products of a 3-cycle and a 4-cycle), namely the 3 in
\eqref{ref} plus
\[\ba{l}
0134625, \ 0135624, \ 0136245, \ 0146235, \ 0235146, \ 0235614, \ 0236145,\\
0245136, \ 0251346, \ 0256134, \ 0346125, \ 0351246, \ 0356124, \ 0361245.
\ea\]

In anticipation of Section \ref{sec:weyl}, consider the action of
$\Si_7$ on the set of $2640$ regular 7-tournaments.  Let
$|T_\si|$ denote the number of regular 7-tournaments fixed
by $\si\in\Si_7$, see after Lemma \ref{lemma:oddlift}. We shall
explain there that there are no fixed tournaments unless $\si$ is the
identity, a 7-cycle or two 3-cycles. In these three cases, $|T_\si|$
equals $\RT(7)=2640$, $2^3$, and $24$ respectively. Burnside's lemma
therefore yields
\[ |\hbox{orbits of $\Si_7$}|=\frac1{7!}\!\!\suml_{\si\in\Si_7}|T_\si|=
\frac1{7!}\big[2640\ti1 + 6!\ti8+
  14\hbox{\small$6\choose3$}\ti24\big]=3,\] confirming that there are
three regular 7-tournaments up to isomorphism.

The respective stabilizers in $\Si_7$ of $\sT_7,\sT_7',\sT_7''$ have
sizes $7,21,3$, with $\sT_7'\cong \sT_7^\QR$ being the most symmetric, with
automorphism group of order 21, accounting for $7!/21=240$ of the
$2640$ labelled regular tournaments. The third isomorphism class is represented by $\sT_7''$, whose automorphism group has order 3, so
that the class contains $1680$ labelled tournaments. The remaining
$720$ constitute the class of $\sT_7$, whose automorphism group is
generated by the 7-fold rotation.

\subsubsection*{Remarks on the case n$=$9}\label{ss9}

As shown in \eqref{3579}, the circulant tournament $\sT_9$ can be
represented by three 9-cycles and three 3-cycles.

It is natural to ask whether one can define a tournament of order $9$
as the union of nine directed 4-cycles. This idea arises from a
geometrically distinguished $4$-form $\phi\in\La^4(\R^9)$ that is
invariant by the action of $\SO(3)\times\SO(3)$ on the irreducible
space $\R^9=\La^2_+\otimes\La^2_-$ \cite{Fino_Nurowski_2015}.
However, an exhaustive search shows that no choice of directed
4-cycles on those nine 4-sets covers each pair exactly once, so this
is likely a red herring.

There are in fact a total of $\RT(9)=3230080$ labelled regular
tournaments of order 9, so each has the same total of balanced
subgraphs. It is known \cite{Chamberland_Herman_2015} that there are
15 isomorphism classes of regular tournaments of order
$9$. Table~\ref{1962} lists balanced subgraphs of sizes up to $7$ (and
also the middle degree $18$). They are ordered by the numbers of
balanced subgraphs of size 7, which can be conveniently memorized as
years lying in four adjacent centuries. (This spurious observation is
on a par with the link to $\2$ in the case $n=7$ above.)

These years
uniquely characterize the isomorphism classes, and $\sT_9$ appears as
the last row, in its own century.

\begin{table}[H]
{\footnotesize\[\begin{array}{r|rrrrrrrr}
q\to&0&3&4&5&6&7&\quad\cdots\quad&18\\\hline  
&1&30&72&180&696&1962&\cdots&377616\\
&1&30&72&180&693&1971&\cdots&377100\\
&1&30&73&178&693&1972&\cdots&377598\\
&1&30&74&176&691&1982&\cdots&377636\\
&1&30&74&176&691&1985&\cdots&377642\\
&1&30&75&174&687&1992&\cdots&377580\\
&1&30&75&174&683&2007&\cdots&377416\\
&1&30&75&174&683&2010&\cdots&377848\\
&1&30&78&168&680&2031&\cdots&377992\\
&1&30&78&168&676&2037&\cdots&377768\\
&1&30&78&168&672&2046&\cdots&377500\\
&1&30&78&168&672&2049&\cdots&377982\\
&1&30&81&162&651&2106&\cdots&375300\\
&1&30&82&160&653&2116&\cdots&379232\\
\sT_9&1&30&90&144&603&2286&\cdots& 392526\\
\end{array}\]}
\vspace{-5pt}
\caption{Some spinor dimensions by degree for $n=9$}\label{1962}
\end{table}

\section{The Weyl group}\label{sec:weyl}

The group of permutations on $n$ letters $\Si_n$ can be embedded in
the normalizer $N(T^n)$ of $T^n$ in $\U(n)$, and identified by a
quotient with the Weyl group $N(T^n)/T^n$.  On $\FF_n=\U(n)/T^n$, fix the metric induced by the scalar product
\begin{equation}
\label{eqn:normalmetric}
\langle X,Y\rangle=-\Tr XY, \quad X,Y\in\lie u(n);
\end{equation}
then $N(T^n)$ acts orthogonally on the
tangent space at $eT^n$, giving an inclusion
$\Si_n\subset\O(2N)$. By construction, the $T^n$-invariant subspace in
any representation of $\O(2N)$ is a representation of $\Si_n$; an
action on the space of invariant spinors can be obtained as
follows. Recall that $\Pin(2N)$ is the subgroup of the Clifford
algebra over $\R^{2N}$ generated by unit vectors $v$; each generator
$v$ acts on $\R^{2N}$ as a reflection around $v^\perp$, which can be
expressed as
\begin{equation}
 \label{eqn:reflection}
 w\mapsto v\cdot w \cdot v.
\end{equation}
More generally, $g\in \Pin(2N)$ acts on $\R^{2N}$ by twisted conjugation
\begin{equation}
 \label{eqn:twisted}
 w\mapsto \theta(g)\cdot w\cdot g^{-1},
\end{equation}
where $\theta$ is the endomorphism of the Clifford algebra extending $\theta(v)=-v$. This gives a homomorphism $\Pin(2N)\to\O(2N)$; we will denote by
$\hat\Si_n$ the preimage of $\Si_n\subset \O(2N)$. Explicitly, extend the definition of the symbols $\alpha_{ij}$ to all pairs $(i,j)$ by declaring
\begin{equation}
 \label{eqn:ijji}
 \alpha_{ji}=-\overline\alpha_{ij};
\end{equation}
a transposition $(i,i+1)$ acts as a reflection on the real plane
$[\Span{\alpha_{i,i+1},\alpha_{i+1,i}}]$, and swaps pairs of real
planes
\[[\![\ol{L_i}\otimes L_k]\!], \quad [\![\ol{L_{i+1}}\otimes L_k]\!], \qquad k\neq i,i+1;\]
the swap of two orthogonal planes is the product of two reflections. The resulting product of $2n-3$ reflections is realized by the product of as many orthogonal unit vectors in $T_e\FF_n$, namely
\[t_j=
 \frac1{\sqrt2}(\alpha_{j,j+1}-\alpha_{j+1,j})\prod_{k\neq j,j+1}
 r_k\,r'_k,
\]
where, for each $k\neq j,j+1$,
\[r_k=\tfrac12(\alpha_{jk}-\alpha_{kj}-\alpha_{j+1,k}+\alpha_{k,j+1}),\qquad
r'_k=\tfrac i2(\alpha_{jk}+\alpha_{kj}-\alpha_{j+1,k}-\alpha_{k,j+1}).\]
Notice that these elements are defined in the complex Clifford algebra, but \eqref{eqn:ijji} forces all $2n-3$ vectors to be real. Consequently $(t_j)^2=(-1)^{(2n-3)(n-1)}=1$, because $n$ is odd.

The group  $\hat \Si_n$ is generated by $t_0,\dotsc, t_{n-2}$ and $z=-1$; $z$ is central of order two, and the following relations hold:
\[z^2=1, \quad t_j^2=1, \quad (t_jt_{j+1})^3=1, \quad t_jt_k=zt_kt_j (j+1<k).\]
The map sending $t_j$ to the transposition $(j, j+1)$ gives an exact sequence
\[1\to \{\pm1\}\to\hat\Si_n\to\Si_n\to 1.\]

\begin{remark}
The image of $t_j$ in $\O(2N)$ reverses orientation, since it compounds $2n-3$ reflections.
\end{remark}

\begin{remark}\label{rem:cliffordcovers}
By construction, $\hat\Si_n$ acts linearly on $\Xi_n$.
Relative to $\Si_n$, $\Xi_n$ is only a projective representation. It is known (see \cite[Theorems 2.9 and 2.12]{Hoffman_Humphreys_1992}) that all projective representations of $\Si_n$ can be lifted to linear representations of one of its Schur covers. For $n\geq 4$, $\Si_n$ has two Schur covers $\tilde\Si_n$ and $\hat\Si_n$, distinguished by whether transpositions lift to elements of order four or two; $\Si_3$ is its own Schur cover, but we will extend the definitions to $n=3$ for a uniform treatment. Both $\tilde\Si_n$ and $\hat\Si_n$ can be realized  inside the complex Clifford algebra of $\R^n$, taking $-1$ and either
\[\frac 1{\sqrt2}(e_j-e_{j+1}), \frac i{\sqrt2}(e_j-e_{j+1})\]
as generators; under twisted conjugation, both act as $\Si_n$ on $\R^n$. A complex spin representation of either group induces a representation of the other by letting e.g. $\frac i{\sqrt2}(e_j-e_{j+1})$ act as $i$ times the action of $\frac 1{\sqrt2}(e_j-e_{j+1})$.

In this description, the intersection of $\tilde\Sigma_n$ and $\hat\Sigma_n$ is the $2:1$ cover of the group of even permutations, and in this intersection two representations in the above correspondence coincide.
\end{remark}

\begin{lemma}\label{lemma:oddlift}
If $\sigma$ has odd order $d$ in $\Si_n$, then it lifts to a unique element $\tilde \sigma$ of odd order in $\hat\Si_n$; the order of $\tilde\sigma$ is $d$.
\end{lemma}
\begin{proof}
There are two elements in $\hat\Si_n$ covering $\sigma$, say $\tilde \sigma$ and $z\tilde \sigma$. Each has necessarily order $d$ or $2d$. Since $d$ is odd,
\[(z\tilde \sigma)^d = z^d\tilde \sigma^d=z \tilde\sigma^d;\]
thus, exactly one has order $d$.
\end{proof}

Every permutation in $\Si_n$ decomposes as the union of cycles; we will say a cycle is odd or even according to whether its order is odd or even. Given a tournament $\sT$, we say a permutation $\sigma$ leaves $\sT$ invariant if it maps arrows to arrows, i.e. if it is an automorphism of the tournament. We denote by $T_\sigma$ the set of $\sigma$-invariant regular tournaments. If $T_\sigma$ is nonempty, $\sigma$  only has odd cycles: indeed, if $\sigma$ has a cycle of order $2k$, then $\sigma^k$ contains a transposition, and therefore cannot leave a tournament invariant.
\begin{remark}
The unique lift with odd order of a permutation with only odd cycles is the product of the odd lifts of each cycle, which commute.
\end{remark}

A permutation $\sigma$ can also be viewed as an automorphism of the undirected complete graph on $n$ nodes; this induces a permutation on the set of undirected edges, call it $\rho(\sigma)$. If $\sigma$ leaves a tournament invariant, $\rho(\sigma)$ can be identified with a permutation on the set of edges of the tournament.

\begin{lemma}\label{lemma:edgecycles}
Let $\sigma$ be a permutation in $\Si_n$ that only has odd cycles. The cycles of $\rho(\sigma)$ are determined  as follows:
\begin{enumerate}
\item two distinct cycles of lengths $a$ and $b$ contribute $\gcd(a,b)$ cycles of length $\operatorname{lcm}(a,b)$, formed by the edges joining them;
\item a cycle of odd length $a$ contributes $(a-1)/2$ cycles of length $a$, formed by the edges joining its nodes.
\end{enumerate}
In particular, $\rho(\sigma)$ only has odd cycles.
\end{lemma}
\begin{proof}
Let $A$ and $B$ be distinct cycles of $\sigma$, of lengths $a$ and $b$. An edge joining a node of $A$ to a node of $B$ is mapped by $\rho(\sigma)$ and each power $\rho(\sigma)^k$ to an edge of the same form; since the endpoints cannot be interchanged, the edge returns to itself when both endpoints do, i.e. when $k$ is a multiple of $\operatorname{lcm}(a,b)$. Hence the $ab$ edges joining $A$ to $B$ split into $\gcd(a,b)$ cycles of $\rho(\sigma)$, each of length $\operatorname{lcm}(a,b)$.

For the edges with both endpoints in a single cycle $A$ of length $a$, identify $A$ with $\Z/a\Z$ so that $\sigma$ acts as $x\mapsto x+1$; an edge is then a pair $\{x,x+d\}$, and the edges
\[\{x,x+d\}, \dotsc, \{x+(a-1),x+(a-1)+d\}\]
are distinct because $a$ is odd.

The final statement follows trivially.
\end{proof}
\begin{proposition}\label{prop:trace}
Let $\sigma$ be a permutation in $\Si_n$.
\begin{itemize}
 \item if $\sigma$ contains even cycles and $\tilde\sigma$ is any lift to $\hat\Si_n$, $\Tr (\tilde \sigma|_{\Xi_n})=0$;
 \item if $\sigma$ only contains odd cycles, then its odd-order lift  satisfies $\Tr (\tilde \sigma|_{\Xi_n})=\abs{T_\sigma}$.
 \end{itemize}
\end{proposition}
\begin{proof}
We compute the trace by evaluating on the basis of pure invariant spinors. If $\psi$ is pure and invariant, then so is $\tilde\sigma(\psi)$, and $\psi$ only gives a contribution to the trace of $\tilde\sigma$ if $\tilde\sigma(\psi)=c\psi$ for some complex scalar $c$; this means that $\tilde\sigma(\psi)$ determines the same almost complex structure as $\psi$, or  equivalently the same regular tournament.  Thus, we only need to consider pairs $(\sigma,\psi)$ such that the regular tournament associated to $\psi$ is in $T_\sigma$. In particular, $\sigma$ only has odd cycles, hence odd order $d$; by Lemma~\ref{lemma:oddlift}, it has a unique odd-order lift $\tilde\sigma$, whose order is also $d$.

Fix a tournament in $T_\sigma$, with associated pure spinor $\psi$, and use the corresponding almost complex structure to identify spinors with $\Lambda^{0,\bullet}$, so that $\psi=1$. Let $W$ be the real span of the products of the $\overline{\alpha_{ij}}$ with $i\to j$ in this tournament, so that $\Lambda^{0,\bullet}=W\oplus iW$; Clifford multiplication by each $\alpha_{ij}$ maps $W$ to itself. Expanding the generators $t_j$, we see that $t_j\cdot W\subset iW$, because $n-2$ is odd. Since $\sigma$ is an even permutation, $\tilde\sigma$ is a product of an even number of generators up to sign, so $\tilde\sigma\cdot W\subset W$ and $c$ is real. Then  $\psi=\tilde\sigma^d(\psi)=c^d\psi$ forces $c=1$, because $d$ is odd; thus, each element of $T_\sigma$ contributes $1$ to the trace.
\end{proof}
We can express $\abs{T_\sigma}$ by a generating function; in order to do so, we introduce the notation $[\,y_1^{e_1}\cdots y_m^{e_m}\,]P$ to denote the coefficient of the monomial $\,y_1^{e_1}\cdots y_m^{e_m}$ in the polynomial $P$.
\begin{lemma}\label{lemma:count}
Let $\sigma$ be a permutation in $\Si_n$ whose cycles have odd lengths $a_1,\dotsc, a_m$, and set $g_{ij}=\gcd(a_i,a_j)$, $L=\sum_i (a_i-1)/2$. Then
\[\abs{T_\sigma}=2^{L}\bigl[y_1^{(n-a_1)/2}\cdots y_m^{(n-a_m)/2}\bigr]\prod_{1\le i<j\le m}\Bigl(y_i^{a_j/g_{ij}}+y_j^{a_i/g_{ij}}\Bigr)^{g_{ij}}.\]
\end{lemma}
\begin{proof}
Write $A_1,\dotsc, A_m$ for the cycles of $\sigma$, viewed as subsets of the set of nodes. Since all the $a_i$ are odd, the proof of Lemma~\ref{lemma:edgecycles} shows that whenever a power of $\sigma$ maps an edge to itself, it fixes both endpoints; hence every $\rho(\sigma)$-orbit of edges can be directed $\sigma$-equivariantly, in exactly two ways, each determined by the direction of any single edge in the orbit. The $\sigma$-invariant tournaments are obtained by directing each orbit independently, and we must count the choices for which the resulting tournament is regular.

By the first item in Lemma~\ref{lemma:edgecycles}, nodes of distinct $A_i$ and $A_j$ are connected by $g_{ij}$ edge orbits, each containing $a_ia_j/g_{ij}$ edges. For each edge orbit, assigning the direction from $A_i$ to $A_j$ contributes $a_j/g_{ij}$ to the out-degree of each node of $A_i$. Thus, the choice corresponds to the generating function $y_i^{a_j/g_{ij}}+y_j^{a_i/g_{ij}}$ in the polynomial ring $\R[y_1,\dotsc, y_m]$. Since there are $g_{ij}$ independent choices for every $i<j$, the overall generating function is
\[
 \prod_{1\le i<j\le m}\Bigl(y_i^{a_j/g_{ij}}+y_j^{a_i/g_{ij}}\Bigr)^{g_{ij}},\]
where the coefficient of $y_1^{e_1}\cdots y_m^{e_m}$ counts choices contributing $e_i$ to the out-degree of nodes of $A_i$.

By the second item of Lemma~\ref{lemma:edgecycles}, the edges connecting nodes of a fixed $A_i$ form $(a_i-1)/2$ orbits,  giving $2^{(a_i-1)/2}$ choices. Each choice amounts to fixing a $\Z/a_i\Z$-invariant tournament on $\Z/a_i\Z$, which is necessarily regular, i.e. the out-degree of each node is $(a_i-1)/2$. The generating function is therefore $2^{(a_i-1)/2} y_i^{(a_i-1)/2}$; taking the product and extracting the coefficient of $y_1^{(n-1)/2}\dotsm y_m^{(n-1)/2}$ gives a formula equivalent to the one in the statement.
\end{proof}

\begin{example}\label{ex:count}
For the $n$-cycle, Lemma~\ref{lemma:count} trivially gives $\abs{T_\sigma}=2^{(n-1)/2}$.

The class $(3^3)$ at $n=9$ exercises the method at full strength: since $g_{ij}=3$ and $a_j/g_{ij}=1$ for all $i,j$, Lemma~\ref{lemma:count} reads
\[\abs{T_\sigma}=2^3\,\bigl[y_1^3y_2^3y_3^3\bigr]\,(y_1+y_2)^3(y_1+y_3)^3(y_2+y_3)^3.\]
Expanding powers into monomials, each monomial of $(y_1+y_2)^3$ only gives a contribution if paired with fixed choices in the expansion of the other two factors, leading to
\[\bigl[y_1^3y_2^3y_3^3\bigr]\,(y_1+y_2)^3(y_1+y_3)^3(y_2+y_3)^3
=\sum_{k=0}^3\binom3k^3=56,\]
so $\abs{T_\sigma}=2^3\cdot56=448$, matching Table~\ref{tab:support}.

Similarly, for $(5^2,3)$ at $n=13$,
\[\abs{T_\sigma}=2^5\,\bigl[y_1^4y_2^4y_3^5\bigr]\,(y_1+y_2)^5(y_1^3+y_3^5)(y_2^3+y_3^5)=
2^5\biggl(\binom51+\binom54\biggr)=320.\]
For $(3,1^2)$ at $n=5$, $T_\sigma$ is empty, because
\[[y_1 y_2^2y_3^2] (y_1+y_2^3)(y_1+y_3^3)(y_2+y_3)=0.\]
\end{example}

\section{The decomposition of \texorpdfstring{$\Xi_n$}{Xi(n)}}\label{sec:decomposition}

Having computed the character of $\hat\Si_n$ on $\Xi_n$, we now
determine the decomposition into irreducible representations. Even
though the character can be determined in terms of $\Si_n$, the
decomposition only contains representations in which $z$ acts as $-1$,
called \emph{spin} representations.

By a theorem of Schur \cite{Schur_1911} (see also \cite[Theorem 8.6]{Hoffman_Humphreys_1992} and \cite{Stembridge_1989}), the irreducible spin representations of $\hat\Si_n$ are parameterized by the \emph{strict} partitions of $n$, i.e. the partitions $\lambda=(\lambda_1,\dotsc,\lambda_\ell)$ with $\lambda_1>\dotsb>\lambda_\ell$. Write $\varepsilon(\lambda)\in\{0,1\}$ for the parity of $n-\ell$. When $\varepsilon(\lambda)=0$, the partition $\lambda$ labels a single spin representation, isomorphic to its own twist by the sign character of $\Si_n$; when $\varepsilon(\lambda)=1$, it labels a pair $\lambda^\pm$ of representations exchanged by the sign twist. In either case the dimension is
\[2^{(n-\ell-\varepsilon(\lambda))/2}\,g_\lambda,\qquad g_\lambda=\frac{n!}{\lambda_1!\dotsm\lambda_\ell!}\prod_{i<j}\frac{\lambda_i-\lambda_j}{\lambda_i+\lambda_j},\]
where $g_\lambda$ counts the standard shifted Young tableaux of shape $\lambda$ \cite[Theorem 10.7]{Hoffman_Humphreys_1992}. By the correspondence of Remark~\ref{rem:cliffordcovers}, this applies to both the Schur covers $\hat\Sigma_n$ and $\tilde\Sigma_n$. By Proposition~\ref{prop:trace}, the character of $\Xi_n$ is zero on lifts of odd permutations;  this implies that $\Xi_n$ is isomorphic to its twist by the sign character, hence two members of a pair $\lambda^\pm$ occur in $\Xi_n$ with the same multiplicity. We accordingly write $\chi_\lambda$ for the character of the representation labelled $\lambda$, or of either member of the pair, and $m_\lambda$ for its multiplicity in $\Xi_n$. Given a partition $\rho$ of $n$, $\Si_n$ contains  $n!/z_\rho$ permutations of type $\rho$, where $z_\rho$ is the size of each centralizer; if $r_j$ denotes the number of times $j$ appears in the partition, $z_\rho =\prod_j j^{r_j}r_j!$. We will write $\abs{T_\rho}$ for the size of $T_\sigma$ with any $\sigma$ of type $\rho$.
\begin{corollary}\label{cor:multiplicity}
For every strict partition $\lambda$ of $n$,
\[
m_\lambda=\sum_{\rho}\frac{\abs{T_\rho}}{z_\rho}\,\overline{\chi_\lambda(\rho)},
\]
where $\rho$ ranges over partitions of $n$ into odd parts and $\chi_\lambda(\rho)$ denotes the value of $\chi_\lambda$ on the odd-order lift of a permutation of cycle type $\rho$.
\end{corollary}
\begin{proof}
The standard projection formula can be written as
\begin{equation}\label{eq:multiplicities}
  m_\lambda \;=\;
  \frac1{\abs{\hat\Si_n}}\sum_{\tilde\sigma\in\hat\Si_n}
  \chi_{\Xi_n}(\tilde\sigma) \overline{\chi_\lambda(\tilde\sigma)}.
\end{equation}
Since $\chi_{\Xi_n}(\tilde\sigma)\overline{\chi_\lambda(\tilde\sigma)}=
  \chi_{\Xi_n}(z\tilde\sigma)\overline{\chi_\lambda(z\tilde\sigma)}$, we have
\[
 m_\lambda \;=\;
  \frac1{\abs{\Si_n}}\sum_{\sigma\in\Si_n}
  \chi_{\Xi_n}(\tilde\sigma)\overline{\chi_\lambda(\tilde\sigma)}.
\]
By Proposition~\ref{prop:trace}, the only $\tilde\sigma$ contributing to the sum are those projecting to permutations that decompose into odd cycles; choosing the odd-order lift gives the formula in the statement.
\end{proof}

The characters $\chi_\lambda$ can be computed using the Morris rule suitably adjusted to odd-order lifts (\cite[Theorem 10.1]{Hoffman_Humphreys_1992}), and are tabulated up to $n=13$ in \cite{Curtis1985ATLASOF}; no distinction between $\hat\Sigma_n$ and $\tilde\Sigma_n$ appears on odd-order lifts, which belong to the common subgroup of Remark~\ref{rem:cliffordcovers}. For $n=3$ and $n=5$ they are small
enough to write out; the values, and the other data entering
Corollary~\ref{cor:multiplicity}, are collected in
Table~\ref{tab:spinchar}.

\begin{table}[t]
  \small
  \centering
 \begin{tabular}[t]{lrrr}
$\lambda$ &   $\chi_\lambda(1^3)$ & $\chi_\lambda(3)$\\
\midrule
$(3)$        & $2$ & $-1$\\
$(2,1)^\pm$  & $1$ & $1$\\
\end{tabular}
\hskip50pt
\begin{tabular}[t]{lrrrr}
$\lambda$ &  $\chi_\lambda(1^5)$ & $\chi_\lambda(3,1^2)$ & $\chi_\lambda(5)$\\
\midrule
$(5)$        & $4$ & $-2$ & $-1$\\
$(4,1)^\pm$  & $6$ & $0$  & $1$\\
$(3,2)^\pm$  & $4$ & $1$  & $-1$\\
\end{tabular}
 \caption{\parshape 1 70pt 320pt
  The irreducible spin characters of $\hat\Si_3$ and
  $\hat\Si_5$, indexed by the strict partitions $\lambda$, evaluated
  on the odd-order lifts of the classes $\rho$ whose parts are all
  odd.}
\label{tab:spinchar}
\end{table}

\begin{example}
For $n=3$,  both regular tournaments
of order three are invariant under a $3$-cycle, so
$\abs{T_{(3)}}=\abs{T_{(1^3)}}=2$, and Corollary~\ref{cor:multiplicity} gives
\[m_{(3)}=\tfrac26\cdot2+\tfrac23(-1)=0,\qquad
m_{(2,1)^\pm}=\tfrac26+\tfrac23=1.\]
For $n=5$ only $\rho=(5)$ and $\rho=(1^5)$ have $\abs{T_\rho}\ne0$ (see Example~\ref{ex:count});
the values of Table~\ref{tab:spinchar} give
\[m_{(3,2)}=m_{(5)}=\frac{24}{120} \cdot 4 + \frac45(-1)=0,\quad
m_{(4,1)}=\frac{24}{120}\cdot6+\frac45\cdot1=2.
\]
For $7\le n\le13$, every spin representation occurs (see Table~\ref{tab:decomp} below).
\end{example}

This leads to the following tables; the traces, obtained from
Lemma~\ref{lemma:count} as in Example~\ref{ex:count}, and the
multiplicities, obtained from Corollary~\ref{cor:multiplicity}, have
been verified by computer for all $n\le13$, in agreement with
Proposition~\ref{prop:trace}:

\begin{table}[H]
  \centering
  \normalsize
  \setlength{\tabcolsep}{5pt}
  \renewcommand{\arraystretch}{1.4}
\begin{tabular}[t]{r r r}
\toprule
$n$ & Cycle type & $\Tr(\tilde\sigma\mid\Xi_n)$ \\
\midrule
\multicolumn{3}{l}{$n=5$, $\RT(5)=24$;} \\
& $(5)$ & $4$ \\
& $(1^5)$ & $24$ \\
\midrule
\multicolumn{3}{l}{$n=7$, $\RT(7)=2640$;} \\
& $(7)$ & $8$ \\
& $(3^2,1)$ & $24$ \\
& $(1^7)$ & $2640$ \\
\midrule
\multicolumn{3}{l}{$n=9$, $\RT(9)=3{,}230{,}080$;} \\
& $(9)$ & $16$ \\
& $(3^3)$ & $448$ \\
& $(3^2,1^3)$ & $160$ \\
& $(3,1^6)$ & $160$ \\
& $(1^9)$ & $3{,}230{,}080$ \\
\bottomrule
\end{tabular}
\quad
\begin{tabular}[t]{r r r}
\toprule
$n$ & Cycle type & $\Tr(\tilde\sigma\mid\Xi_n)$ \\
\midrule
\multicolumn{3}{l}{$n=11$, $\RT(11)=48{,}251{,}508{,}480$;} \\
& $(11)$ & $32$ \\
& $(5^2,1)$ & $320$ \\
& $(3^2,1^5)$ & $6720$ \\
& $(3,1^8)$ & $215{,}040$ \\
& $(1^{11})$ & $48{,}251{,}508{,}480$ \\
\midrule
\multicolumn{3}{l}{$n=13$, $\RT(13)=9{,}307{,}700{,}611{,}292{,}160$;} \\
& $(13)$ & $64$ \\
& $(5^2,3)$ & $320$ \\
& $(5^2,1^3)$ & $2240$ \\
& $(3^4,1)$ & $470{,}016$ \\
& $(3^2,1^7)$ & $2{,}943{,}360$ \\
& $(3,1^{10})$ & $3{,}053{,}998{,}080$ \\
& $(1^{13})$ & $9{,}307{,}700{,}611{,}292{,}160$ \\
\bottomrule
\end{tabular}\medskip
\caption{\parshape 1 70pt 330pt
  Contributing cycle types for $n\le 13$:\ all types with
  $\Tr(\tilde\sigma\mid\Xi_n)\ne 0$ for the odd-order lifts
  $\tilde\sigma$. Other types give zero.}
\label{tab:support}
\end{table}

\begin{table}[H]
\centering
\small
\setlength{\tabcolsep}{5pt}
\renewcommand{\arraystretch}{1.05}
\begin{tabular}[t]{l r r}
\toprule
$\lambda$ & dimension & $m_\lambda$ \\
\midrule
\multicolumn{3}{l}{$n=3$;}\\
$(2,1)^\pm$ & $1$ & $1$ \\
\midrule
\multicolumn{3}{l}{$n=5$;}\\
$(4,1)^\pm$ & $6$ & $2$ \\
\midrule
\multicolumn{3}{l}{$n=7$;}\\
$(7)$ & $8$ & $8$ \\
$(6,1)^\pm$ & $20$ & $8$ \\
$(5,2)^\pm$ & $36$ & $20$ \\
$(4,3)^\pm$ & $20$ & $12$ \\
$(4,2,1)$ & $28$ & $12$ \\
\midrule
\multicolumn{3}{l}{$n=9$;}\\
$(9)$ & $16$ & $144$ \\
$(8,1)^\pm$ & $56$ & $504$ \\
$(7,2)^\pm$ & $160$ & $1416$ \\
$(6,3)^\pm$ & $224$ & $1992$ \\
$(6,2,1)$ & $240$ & $2144$ \\
$(5,4)^\pm$ & $112$ & $1016$ \\
$(5,3,1)$ & $336$ & $2976$ \\
$(4,3,2)$ & $96$ & $872$ \\
\midrule
\multicolumn{3}{l}{$n=11$;}\\
$(11)$ & $32$ & $38{,}688$ \\
$(10,1)^\pm$ & $144$ & $174{,}016$ \\
$(9,2)^\pm$ & $560$ & $676{,}752$ \\
$(8,3)^\pm$ & $1200$ & $1{,}450{,}352$ \\
$(8,2,1)$ & $1232$ & $1{,}489{,}008$ \\
$(7,4)^\pm$ & $1440$ & $1{,}740{,}624$ \\
$(7,3,1)$ & $3168$ & $3{,}829{,}312$ \\
$(6,5)^\pm$ & $672$ & $812{,}368$ \\
$(6,4,1)$ & $3168$ & $3{,}829{,}680$ \\
$(6,3,2)$ & $2464$ & $2{,}978{,}688$ \\
$(5,4,2)$ & $1760$ & $2{,}127{,}744$ \\
$(5,3,2,1)^\pm$ & $528$ & $638{,}320$ \\
\bottomrule
\end{tabular}
\quad
\begin{tabular}[t]{l r r}
\toprule
$\lambda$ & dimension & $m_\lambda$ \\
\midrule
\multicolumn{3}{l}{$n=13$;}\\
$(13)$ & $64$ & $95{,}655{,}104$ \\
$(12,1)^\pm$ & $352$ & $526{,}109{,}088$ \\
$(11,2)^\pm$ & $1728$ & $2{,}582{,}770{,}624$ \\
$(10,3)^\pm$ & $4928$ & $7{,}365{,}802{,}016$ \\
$(10,2,1)$ & $4992$ & $7{,}461{,}463{,}232$ \\
$(9,4)^\pm$ & $8800$ & $13{,}153{,}378{,}880$ \\
$(9,3,1)$ & $18{,}304$ & $27{,}359{,}081{,}664$ \\
$(8,5)^\pm$ & $9504$ & $14{,}205{,}773{,}440$ \\
$(8,4,1)$ & $32{,}032$ & $47{,}878{,}866{,}112$ \\
$(8,3,2)$ & $22{,}464$ & $33{,}577{,}437{,}696$ \\
$(7,6)^\pm$ & $4224$ & $6{,}313{,}707{,}456$ \\
$(7,5,1)$ & $27{,}456$ & $41{,}039{,}271{,}424$ \\
$(7,4,2)$ & $41{,}600$ & $62{,}180{,}898{,}272$ \\
$(7,3,2,1)^\pm$ & $9152$ & $13{,}679{,}819{,}264$ \\
$(6,5,2)$ & $22{,}464$ & $33{,}577{,}808{,}288$ \\
$(6,4,3)$ & $18{,}304$ & $27{,}359{,}755{,}936$ \\
$(6,4,2,1)^\pm$ & $13{,}728$ & $20{,}519{,}840{,}880$ \\
$(5,4,3,1)^\pm$ & $4576$ & $6{,}839{,}971{,}360$ \\
\bottomrule
\end{tabular}\medskip
\caption{\parshape 1 70pt 350pt
Decomposition of $\Xi_n$ into spin representations for $n\le 13$:
dimension and multiplicity of every representation that occurs. A row
$\lambda^\pm$ refers to both members of an associate pair, each with
the stated dimension and multiplicity.\label{tab:decomp}}
\end{table}

Every multiplicity appearing in Table~\ref{tab:decomp} is even, except
for $n=3$. We will see in Corollary~\ref{cor:mod4} that this holds
generally for $n\equiv 1 \pmod 4$.

\begin{proposition}
Let $\lambda$ be a strict partition of $n$ with $\ell$ parts, and write
$k=n-\ell$. Suppose that, modulo eight,
\[n\equiv 1,3\ \text{ and }\ k\equiv 4,5,6,\qquad\text{or}\qquad
n\equiv 5,7\ \text{ and }\ k\equiv 0,1,2.\]
Then $m_\lambda$ is even.
\end{proposition}

\begin{proof}
It is a general fact that the multiplicity of a representation $V$ in
$W$ is even if one of $V,W$ is real and the other quaternionic,
because $V^*\otimes W$ is quaternionic.

By Schur \cite[Satz III]{Schur_1927} the representation
labelled $\lambda$ is of real type when $k\equiv 0,1,2\pmod 8$, of complex
type when $k\equiv 3,7$, and of quaternionic type when $k\equiv 4,5,6$.

On the other hand, the spinor representation $\Delta$ is real when
$2N\equiv 0,6\pmod 8$ and quaternionic when $2N\equiv 2,4$; the
actions of $T^n\subset\Spin(2N)$ and $\hat\Sigma_n\subset \Pin(2N)$
preserve the structure, hence $\Xi_n$ is a
$\hat\Sigma_n$-representation of real type when $n\equiv 1,3\pmod 8$
and of quaternionic type when $n\equiv 5,7$. The statement follows.
\end{proof}

\section{The Dirac operator}\label{sec:dirac}
Extend the scalar product \eqref{eqn:normalmetric} bilinearly to $\gl(n,\C)$; the generators $\alpha_{hj}$ satisfy
\[\langle \alpha_{hj},\alpha_{kl}\rangle=-\delta_{jk}\delta_{hl}.\]
The conjugation of $\gl(n,\C)$ relative to the real form $\lie u(n)$ is $X\mapsto -X^*$, so that $\overline{\alpha_{hj}}=-\alpha_{jh}$ and  $\langle \alpha_{hj},\overline{\alpha_{hj}}\rangle = 1$. Thus,  each $\alpha_{hj}$ is a unit vector for the Hermitian metric, and the vector dual to $\alpha_{hj}$ under
$\langle,\rangle$ is $-\alpha_{jh}$.

Write $\lie u(n)=\lie t^n\op \lie m$, where $\lie t^n$ consists of the diagonal matrices and $\lie m\otimes\C$ is spanned by the $\alpha_{hj}$, and denote by subscript $\lie m$ the projection. For $X,Y,Z$ in $\lie m$, define
\[B(X,Y,Z)=\langle [X,Y],Z\rangle.\]
Up to a negative factor and a restriction, $B$ is the Cartan $3$-form. If $i,j,l$ are distinct, then $[\alpha_{ij},\alpha_{jl}]=-\alpha_{il}$, so that $B(\alpha_{ij},\alpha_{jl},\alpha_{li})=1$; up to permutations, these are the only nonzero values of $B$ on triples of elementary matrices. Since the vectors $\alpha_{ij}$, $\alpha_{jl}$, $\alpha_{li}$ are mutually orthogonal, $B$ corresponds to the element
\[B=\frac13\sum_{i,j,l}\alpha_{ij}\cdot\alpha_{jl}\cdot\alpha_{li}\]
of the Clifford algebra, where the sum is over all triples of distinct indices.

The metric $g=\langle,\rangle$ is normal homogeneous, hence naturally
reductive, and the invariant spinors are the constant spinors of
\cite{Agricola_2003}.  In the
notation of that paper, the connections
\[\nabla^t_XY=\nabla^0_XY+t[X,Y]_{\lie m}\]
form a family joining the canonical connection $\nabla^0$ to the Levi-Civita connection $\nabla^{1/2}$, the cubic element is $H=\frac32B$, and $D^t$ denotes the Dirac operator of $\nabla^t$. On invariant spinors \cite[Theorem 4.2]{Agricola_2003} gives
\begin{equation}
 \label{eqn:derivativeofspinors}
 \nabla^t_X\psi = \frac t3 (X\hook H)\cdot\psi = \frac t2 (X\hook B)\cdot\psi,
\end{equation}
hence,
\[D^t\psi=\frac t2\sum_i e_i\cdot(e_i\hook B)\cdot\psi=\frac{3t}2 B\cdot\psi.\]
Since all $D^t$ are proportional, we shall fix the Levi-Civita Dirac operator $D=D^{1/2}$. Consistently with \cite[Theorem 4.2]{Agricola_2003}, we show that its square acts as a multiple of the identity.
\begin{proposition}
\label{prop:eigenspinors}
$\Xi_n$ splits as the sum of two eigenspaces for $D$ with eigenvalue
$\pm\frac14\sqrt{3n(n^2-1)}$, each of
dimension $\tfrac12\RT(n)$.
\end{proposition}
\begin{proof}
Since $D$ interchanges chirality, $\Tr D=0$; it is therefore sufficient to show that
\begin{equation}
\label{eqn:squareD}
D^2=\frac{3n(n^2-1)}{16}\,\id.
\end{equation}
Since $\Xi_n$ is generated by pure spinors, it suffices to compute $D^2\psi$ for $\psi$ pure. We can use the complex structure defined by $\psi$ itself, so that $\psi=1$ as a form. Then
\[D\psi=\frac1{\sqrt2}\sum_{i,j,l}\alpha_{ij}\we \alpha_{jl}\we \alpha_{li},\]
where the sum is over $i,j,l$ such that $\alpha_{ij}$, $\alpha_{jl}$, $\alpha_{li}$ have type $(0,1)$; equivalently, over cycles $jil$ in the tournament $\sT$ defining the fixed complex structure.

The operator $D$ has odd degree, and being cubic, $D^2$ can only change the degree by $0,2,4,6$ (only raising occurs, since $\psi=1$); however, $D^2\psi$ cannot have components of degree two, because there are no balanced subgraphs with two edges. Moreover  $D^2\psi$ has no component of degree six, because three-forms anticommute. Degree four components in $D^2\psi$ arise from terms
\[\alpha_{ik}\we\alpha_{kj}\we \alpha_{ji}\hook (\alpha_{ij}\we \alpha_{jl}\we \alpha_{li})=-\alpha_{ik}\we\alpha_{kj}\we \alpha_{jl}\we \alpha_{li} \]
where $jil$ and $kilj$ are directed cycles in $\sT$. In other words,
this degree $4$ spinor is obtained by joining a directed $3$-cycle in
$\sT$ and one that has an inverted arrow compared to the tournament,
with the two differently oriented arrows `canceling out'.

Now whenever  $jil$ and $kilj$ are directed cycles,  $\sT$ contains
either $kil$ or $klj$, and we obtain corresponding terms
\begin{align*}
\alpha_{kj}\we\alpha_{jl}\we \alpha_{lk}\hook (\alpha_{ik}\we \alpha_{kl}\we \alpha_{li})&=\alpha_{kj}\we\alpha_{jl}\we \alpha_{ik}\we  \alpha_{li},\\
\alpha_{li}\we\alpha_{ik}\we \alpha_{kl}\hook (\alpha_{lk}\we \alpha_{kj}\we \alpha_{jl})&=-\alpha_{li}\we\alpha_{ik}\we \alpha_{kj}\we \alpha_{jl}.
\end{align*}
This amounts to dividing a square into two triangles by choosing a diagonal; the diagonal that is chosen determines the resulting sign, so the overall sum is zero.

Thus, the only surviving terms come from contracting a cycle with its
inverse, for example
\[\alpha_{ji}\cdot\alpha_{il}\cdot\alpha_{lj}\cdot\alpha_{ij}\cdot \alpha_{jl}\cdot\alpha_{li}\cdot 1=(-2\sqrt2 )(-2\sqrt2)=8;
\]
each $3$-cycle determines $9$ summands, because indices can be reordered. Since a regular tournament has $n(n^2-1)/24$ oriented $3$-cycles,
\[D^2\psi = \frac1{16} 72 \abs{\{\text{oriented $3$-cycles}\}} =\frac92\frac{n(n^2-1)}{24}=\frac{3n(n^2-1)}{16}.\]
\end{proof}

Equation~\eqref{eqn:squareD} can also be obtained from the Lichnerowicz
formula $D^2=\nabla^*\nabla+\frac14\scal$. Indeed, a direct computation from \eqref{eqn:derivativeofspinors} gives $\nabla^*\nabla\psi=\frac{n(n-1)^2}{16}\psi$ on $\Xi_n$. On the other hand, by \cite[Proposition~7.87]{Besse_1987} the curvature is
\[\langle R(X,Y)X,Y\rangle=\abs{[X,Y]_{\lie t^n}}^2+\frac14\abs{[X,Y]_{\lie m}}^2, \quad X,Y\in\lie m.\]
Summing over an orthonormal basis of $\lie m$ gives
\[\mathrm{scal}=\frac{n(n-1)(n+2)}2;\]
therefore, the
Lichnerowicz formula gives back
\[D^2\psi=\nabla^*\nabla\psi+\frac14\mathrm{scal}\,\psi=
\Bigl(\frac{n(n-1)^2}{16}+\frac{n(n-1)(n+2)}8\Bigr)\psi=\frac{3n(n^2-1)}{16}\psi.\]

Our next result relates the Dirac operator to the $\hat\Si_n$-decomposition. Clifford multiplication is equivariant under twisted conjugation, and \eqref{eqn:twisted} implies
\[\theta(\tilde \sigma)\cdot v \cdot \tilde\sigma^{-1} = \sigma(v).\]
Since $\theta(\tilde\sigma)=\sign(\sigma)\tilde\sigma$, we obtain
\[\tilde\sigma\cdot v\cdot\tilde\sigma^{-1}=\sign(\sigma)(\sigma(v));\]
the invariance of $B$ (which has degree three) now reads
\[\tilde\sigma \cdot B\cdot \tilde\sigma^{-1}=\sign(\sigma)B.\]
Therefore, the Dirac operator is not $\hat\Si_n$-equivariant. To obtain an honest equivariant operator, we consider
\[A\psi=(-1)^q\frac4{\sqrt{3n(n^2-1)}}D\psi, \quad \psi\in\Xi_{n,q}.\]
\begin{corollary}\label{cor:doubling}
The operator $A\colon\Xi_n\to\Xi_n$ commutes with the action of
$\hat\Si_n$ and satisfies $A^2=-\id$. Its eigenspaces $P^+$ and
$P^-$ are subrepresentations, and
\[\Xi_n\cong P^+\op(P^+\otimes\sign).\]
In particular,
$m_\lambda$ is even for every strict partition $\lambda$ of $n$ with
$\varepsilon(\lambda)=0$.
\end{corollary}
\begin{proof}
Equivariance holds by construction, and $A^2=-\id$ follows from Proposition~\ref{prop:eigenspinors}.

Let $\tau$ be the involution of $\Xi_n$ acting as $(-1)^q$ on $\Xi_{n,q}$. Then
\[\tau(\tilde\sigma \cdot \psi)=\sign(\sigma)\tilde\sigma\cdot\tau(\psi).\]
Thus, $\tau$ maps every irreducible $\hat\Si_n$-representation $W$ to $\sign\otimes W$. On the other hand, $\tau$ anticommutes with $D$, hence $A$, therefore interchanging $P^+$ and $P^-$.

If $\varepsilon(\lambda)=0$, the representation labelled $\lambda$ is
isomorphic to its own twist by the sign character, so it occurs as many times in $P^+$ as it does
in $P^+\otimes\sign$, forcing $m_\lambda$ to be even.
\end{proof}

\begin{corollary}\label{cor:mod4}
If $n\equiv1\pmod 4$, then $m_\lambda$ is even for every strict partition
$\lambda$ of $n$.
\end{corollary}
\begin{proof}
Consider the real element of the Clifford algebra
\[\omega=\prod_{i<j}\frac1{\sqrt2}(\alpha_{ij}-\alpha_{ji}).\]
Having fixed a tournament $\sT$, examining the action of the individual factors shows that Clifford multiplication by $\omega$ maps $\psi_\sS$ to $\pm\psi_{\sT\setminus\sS}$; in other words, $\omega$ realizes the isomorphism $*$ of Lemma~\ref{Pd}.

This isomorphism can be seen to be $\hat\Sigma_n$-equivariant as follows. Each factor $\alpha_{ij}-\alpha_{ji}$ of $\omega$ commutes with itself and anticommutes with $\alpha_{ij}+\alpha_{ji}$; since each of the remaining $N-1$ factors anticommutes with $\alpha_{ij}\pm\alpha_{ji}$, and $N-1$ is odd due to the assumption $n\equiv1\pmod 4$, we obtain
\begin{equation}
 \label{eqn:omegaalphaij}
\begin{aligned}
\omega\cdot (\alpha_{ij}-\alpha_{ji})&=-(\alpha_{ij}-\alpha_{ji})\cdot\omega,\\
\omega\cdot (\alpha_{ij}+\alpha_{ji})&=(\alpha_{ij}+\alpha_{ji})\cdot\omega.
\end{aligned}
\end{equation}
Each $t_j$ is the product of $1+(n-2)$ elements of the form $\alpha_{ij}-\alpha_{ji}$, which anticommute with $\omega$, and $n-2$ elements that commute; since $1+(n-2)$ is even, $t_j$ commutes with $\omega$. Hence, $\omega$ acts in a $\hat\Sigma_n$-equivariant way.

Moreover, \eqref{eqn:omegaalphaij} implies $\omega\cdot\alpha_{ij}=\alpha_{ji}\cdot\omega$, hence
\[\omega\cdot \alpha_{ij}\cdot\alpha_{jl}\cdot\alpha_{li}
=\alpha_{ji}\cdot\alpha_{lj}\cdot\alpha_{il}\cdot\omega
=-\alpha_{ji}\cdot\alpha_{il}\cdot\alpha_{lj}\cdot\omega.
\]

This implies that $\omega$ anticommutes with $D$. Since $N$ is even,
$(-1)^q$ and $(-1)^{N-q}$ coincide and $\omega$ anticommutes with
$A$. Thus, $\omega$ interchanges the eigenspaces $P^+,P^-$ of
Corollary~\ref{cor:doubling}, and each irreducible representation
occurs with even multiplicity.
\end{proof}

{\vskip10pt\small\noindent{\bf Acknowledgements.} This note can be
  regarded as a sequel to \cite{ArtachoSemmelmann}, albeit restricted
  to the unitary case, and we are grateful to the authors of that
  paper and Marie-Am\'elie Lawn for introducing a challenging
  problem. S.S.\ learnt about it as an examiner of \cite{deObeso}, and
  is grateful to Martin Forde and Konstanze Rietsch for help with some
  initial computations and examples. He presented results from this
  paper at the Clay Mathematics Institute conference Hitchin80.

  D.C. acknowledges Gruppo Nazionale per le Strutture Algebriche,
Geometriche e le loro Applicazioni (GNSAGA) of Istituto Nazionale di Alta
Matematica (INdAM), the MIUR Excellence Department Project CUP I57G22000700001 awarded to the Department of Mathematics, University of Pisa and the ChatGPT for Academic Researchers programme.
  }

{\vskip5pt\small\noindent{\bf Generative AI disclosure.} During the
    development of this work, the authors used both Claude and ChatGPT
    to assist with exploratory computations, possible proof
    directions, and editorial polishing. The authors wrote the
    mathematical proofs and computations, and take full responsibility
    for the contents of the paper.}

\bibliographystyle{acm} 
\bibliography{isp}

\end{document}